\documentclass[a4paper,landscape]{article}

\usepackage{a4} 
\usepackage{amsfonts}
\usepackage{amssymb}
\usepackage{amsmath}
\usepackage{amscd}
\usepackage{amsthm}
\usepackage{epsfig}
\usepackage[all]{xy}
\usepackage{color}
\usepackage{bbm}
\usepackage{graphicx}
\usepackage{psfrag}
\usepackage{hyperref}
\usepackage{enumerate}

\newcommand{\R}{\mathbb{R}}

\newcommand{\Z}{\mathbb{Z}}
\newcommand{\Q}{\mathbb{Q}}

\newcommand{\bprf}{\begin{proof}[Proof]}
\newcommand{\eprf}{\end{proof}}

\newcommand{\Ka}{{\sf K}}
\newcommand{\La}{{\sf L}}

\newcommand{\Fe}{{\cal F}}
\newcommand{\Ge}{{\cal G}}

\newcommand{\ra}{\rightarrow}

\newcommand{\conv}{{\rm conv}}
\newcommand{\cone}{{\rm cone}}

\newcommand{\conn}{{\rm conn}}

\graphicspath{{pictures/}}   

\newtheoremstyle{plainNoItalics}{}{}{\normalfont}{}{\bfseries}{.}{ }{}

\newtheorem{theorem}{Theorem}
\newtheorem{proposition}[theorem]{Proposition}

\newtheorem{lemma}[theorem]{Lemma}

\newtheorem*{theorem*}{Theorem}
\newtheorem*{proposition*}{Proposition}
\newtheorem*{lemma*}{Lemma}
\newtheorem*{corollary*}{Corollary}

\theoremstyle{plainNoItalics}

\newtheorem*{definition*}{Definition}

\newtheorem*{remark*}{Remark}

\newtheorem*{observation*}{Observation}
\newtheorem*{example*}{Example}
\newtheorem*{problem*}{Problem}

\begin{document}
\title{Tverberg's theorem with constraints}
\author{{\Large Stephan Hell}}
\date{Institut f\"ur Mathematik, MA 6--2,
TU Berlin,\\ D--10623 Berlin, Germany,
hell@math.tu-berlin.de}

\maketitle

\begin{abstract}
  The topological Tverberg theorem claims that for any continuous map 
  of the $(q-1)(d+1)$-simplex $\sigma^{(d+1)(q-1)}$ to $\R^d$ there 
  are $q$ disjoint faces of $\sigma^{(d+1)(q-1)}$
  such that their images have a non-empty intersection. This has been
  proved for affine maps, and if $q$ is a prime power, but not in
  general.\\
  We extend the topological Tverberg theorem in the following way:
  Pairs of vertices are forced to end up in different faces.
  This leads to the concept of constraint graphs. In Tverberg's
  theorem with constraints, we come up with a list of constraints
  graphs for the topological Tverberg theorem.\\
  The proof is based on
  connectivity results of chessboard-type complexes. Moreover, Tverberg's
  theorem with constraints implies new lower bounds for the number of
  Tverberg partitions. As a consequence, we prove Sierksma's 
  conjecture for $d=2$, and $q=3$.
\end{abstract}
\begin{section}{Introduction}\label{sec-intro}
  Helge Tverberg showed in 1966 that any $(d+1)(q-1)+1$ points in $\R^d$ can
  be partitioned into $q$~subsets such that their convex hulls have a
  non-empty intersection. This has been generalized to the
  following statement by B\'ar\'any et
  al.~\cite{bss81:_tverb} for primes~$q$, and by 
  \"Ozaydin~\cite{oezaydin87:_equiv} and 
  Volovikov~\cite{volovikov96:_tverb_theor}  
  for prime powers $q$, using the equivariant
  method from topological combinatorics. The general case for arbitrary
  $q$ is open.
\begin{theorem}\label{thm-ttt}
  Let $q\geq 2$ be a prime power, $d\geq 1$.  For every continuous map
  \linebreak $f:\|\sigma^{(d+1)(q-1)}\|\rightarrow\R^d$ there
  are $q$ disjoint faces $F_1,F_2,\ldots,F_q$ in the standard
  $(d+1)(q-1)$-simplex $\sigma^{(d+1)(q-1)}$ such that their images
  under $f$ have a non-empty intersection.
\end{theorem}
  The special case for affine maps $f$ is equivalent to the original
  statement of Tverberg.
  A partition $F_1,F_2,\ldots,F_q$ as above is a {\it Tverberg
  partition}. A point in the non-empty intersection is a {\it Tverberg 
  point}. In 2005, Sch\"oneborn and Ziegler
\cite[Theorem 5.8]{schoeneborn05:_topol_tverb_theor} showed that for
primes $p$ every continuous map \linebreak
$f:\|\sigma^{3p-3}\|\rightarrow\R^2$ has a Tverberg partition
subject to the following type of constraints: Certain pairs of points
end up in different partition sets. In other words, there is a Tverberg
partition that does not use the edge connecting this pair of points.

To formalize this, let $G$ be a subgraph of the $1$-skeleton of
$\sigma^{(d+1)(q-1)}$, and $f:\sigma^{(d+1)(q-1)}\ra\R^d$ be a
continuous map. Let $E(G)$ be the set of edges of $G$.  A Tverberg
partition $F_1,F_2,\ldots F_q\subset \sigma^{(d+1)(q-1)}$ of
$f$ is a {\it Tverberg partition of $f$ not using any edge of $G$} if
\[|F_i\cap e|\leq 1\text{ for all }i\in[q]\text{ and all edges }e\in E(G).\]

Their proof can easily be carried over to arbitrary dimension $d\geq
1$, and to prime powers $q$ so that one obtains the following statement.
A {\it matching} on a graph $G$ is a set of edges of $G$ such that no two 
of them share a vertex in common.

\begin{theorem}\label{thm-tv-ohne-matching}Let $q>2$
  be a prime power, and $M$ a matching on the graph of
  $\sigma^{(d+1)(q-1)}$. Then every continuous map 
  $f:\|\sigma^{(d+1)(q-1)}\|\ra\R^d$ has a Tverberg partition
  $F_1,F_2,\ldots,F_q$ not using any edge from $M$.
\end{theorem}
Sch\"oneborn and Ziegler use the more general concept of
winding partitions. For the sake of simplicity, we do not
use this setting. However, all results in this paper also
hold for winding partitions.
 
Theorem~\ref{thm-tv-ohne-matching} was an important step for better
understanding of Tverberg partitions: One can force pairs of points to be
in different partition sets of a Tverberg partition.  Choose disjoint
pairs of vertices of $\sigma^{(d+1)(q-1)}$, then this choice
corresponds to a matching $M$ in the $1$-skeleton of
$\sigma^{(d+1)(q-1)}$. For any map $f$, the endpoints of any edge in
$M$ end up in different partition sets due to 
Theorem~\ref{thm-tv-ohne-matching}. 

We extend their result to a wider class of graphs based on the
following approach.  

\begin{definition*}A {\it constraint graph $C$} in
$\sigma^{(d+1)(q-1)}$ is a subgraph of the graph of $\sigma^{(d+1)(q-1)}$
such that every continuous map
$f:\|\sigma^{(d+1)(q-1)}\|\ra\R^d$ has a Tverberg partition of
disjoint faces not using any edge from $C$.
\end{definition*}

Theorem~\ref{thm-tv-ohne-matching} implies that any matching in
$\sigma^{(d+1)(q-1)}$ is a constraint graph for prime powers $q$.
Sch\"oneborn and Ziegler~\cite{schoeneborn05:_topol_tverb_theor} also
come up with an example showing that the bipartite graph $K_{1,q-1}$
is not a constraint graph for arbitrary $q$.

The alternating drawing of $K_{3q-2}$ is shown in Figure
\ref{fig-k10-ohne-3} for $q=4$. If one deletes the first $q-1$ edges
incident to the right-most vertex, then one can check that there is no
Tverberg partition. In Figure~\ref{fig-k10-ohne-3}, the deleted edges
are drawn in broken lines.  Numbering the vertices from right to left
with the natural numbers in $[3q-2]$, the edges of the form
$(1,3q-2-2i)$, for $0\leq i\leq q-2$, are deleted.
\begin{figure}[h]
  \centering \includegraphics{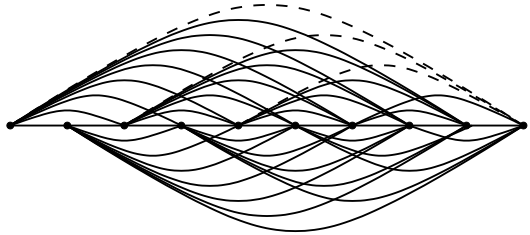} \caption{$K_{10}$
  minus three edges with no winding partition.}
  \label{fig-k10-ohne-3}
\end{figure}

The following theorem generalizes both Theorems~\ref{thm-ttt} and
\ref{thm-tv-ohne-matching}. Moreover, it implies that $K_{1,q-1}$ is a
minimal example for prime powers $q$: All subgraphs of $K_{1,q-1}$ 
are constraint graphs.
\begin{theorem}\label{thm-constraint-graphs}
  Let $q>2$ be a prime power. Then the following subgraphs of
  $\sigma^{(d+1)(q-1)}$ are constraint graphs:
\begin{enumerate}[\rm i)]
\item\label{comp-const}Complete graphs $K_{l}$ on $l$ vertices for $2l<q+2$,
\item\label{bip-const}complete bipartite graphs $K_{1,l}$ for $l<q-1$,
\item\label{path-const}paths $P_l$ on $l+1$ vertices for $l\leq(d+1)(q-1)$ and $q>3$,
\item\label{cycle-const}cycles $C_l$ on $l$ vertices for
  $l\leq(d+1)(q-1)+1$ and $q>4$,
\item\label{union-const}and arbitrary disjoint unions of graphs from
  {\rm (\ref{comp-const})--(\ref{cycle-const})}.
\end{enumerate}
\end{theorem}
The family of constraint graphs is closed under taking subgraphs. It
is thus a monotone graph property.  Theorem~\ref{thm-constraint-graphs} serves us below to estimate the number of
Tverberg points in the prime power case. It is easy to see that $K_2$ is
not a constraint graph for $q=2$.

Figure~\ref{fig-k13-constraint} shows an example of a configuration of
$13$ points in the plane together with a constraint graph. Theorem
\ref{thm-constraint-graphs} implies that there is a Tverberg partition
into $5$ blocks that does not use any of the broken edges. In Figure
\ref{fig-k13-constraint}, there is for example the Tverberg partition
$\{6,10\}$, $\{9,11\}$, $\{0,2,8\}$, $\{1,5,12\}$, $\{3,4,7\}$ that
does not use any of the broken edges. 

The constraint graph $K_l$ guarantees that all $l$ points end up in
$l$ pairwise disjoint partition sets. The constraint graph $K_{1,l}$
forces that the singular point in one shore of $K_{1,l}$ ends up in
a different partition set than all $l$ points of the other shore. \\
\begin{figure}[h]
  \centering
  \includegraphics{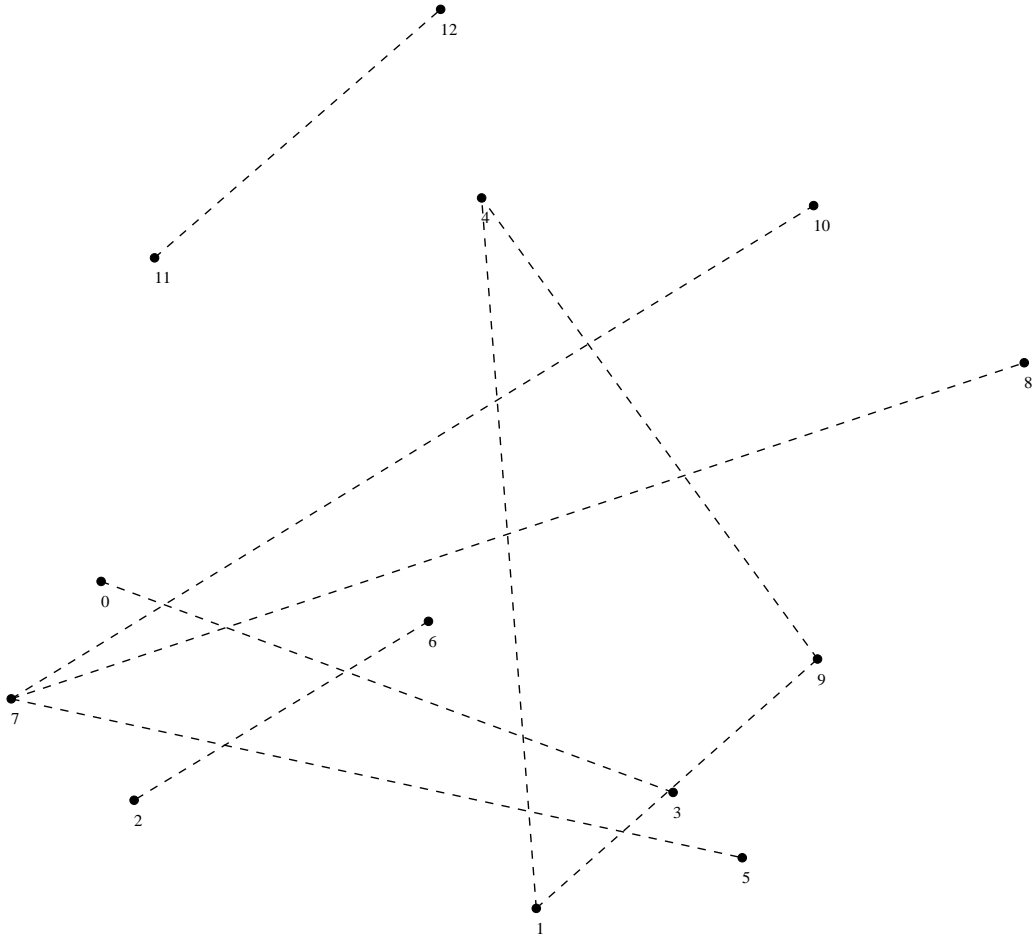}
  \caption{A planar configuration together with a constraint graph for
    $q=5$.}
  \label{fig-k13-constraint}
\end{figure}

{\bf On the number of Tverberg partitions.} Tverberg's theorem
establishes the existence of at least one Tverberg partition. 
 Vu\'ci\'c and \v{Z}ivaljevi\'c~\cite{vz93:_notes_sierk}, and 
Hell~\cite{hell07:_tverb} showed that there is at least 
\[
\frac{1}{(q-1)!}\cdot\left(\frac{q}{r+1}\right)^{\lceil\frac{(d+1)(q-1)}{2}\rceil}\]
many Tverberg partitions if $q=p^r$ is a prime power.

Recently, Hell~\cite{hell07:_birch} showed a lower bound in the original
affine setting of Tverberg which holds for arbitrary~$q$.
\begin{theorem}\label{thm-lower-aff-arb}Let $X$ be a set of
  $(d+1)(q-1)+1$ points in general position in $\R^d$, $d\geq 1$.
  Then the following properties hold for the number $T(X)$ of Tverberg
  partitions: \begin{enumerate}[\rm i)]
  \item $T(X)$ is even for $q>d+1$.  
  \item\label{item-tp-lower} $T(X)\geq (q-d)!$ 
  \end{enumerate}
\end{theorem}
Sierksma conjectured in 1979 that the number of Tverberg partitions is
at least $((q-1)!)^d$. This conjecture is unsettled, except for the
trivial cases $q=2$, or $d=1$.  Using Theorem~\ref{thm-constraint-graphs} 
on Tverberg partitions with constraints  
we can improve the lower bound for the affine setting 
of Theorem~\ref{thm-lower-aff-arb} in
the prime power case.
\begin{theorem}\label{thm-low-aff-improv}Let $d\geq 2$, 
and $q>2$ be a prime power. Then there is an integer
constant $c_{d,q}\geq 2$ such that every set $X$ of $(d+1)(q-1)+1$ points
in general position in $\R^d$ has at least
\[\min\{(q-1)!,\,c_{d,q}(q-d)!\} \]
many Tverberg partitions. Moreover, the constant $c_{d,q}$ is
monotonely increasing in $q$, and $c_{2,3}=4$.
\end{theorem}
This settles Sierksma's conjecture for a wide class of
planar sets for $q=3$. Using some more effort, we entirely
establish Sierksma's conjecture for $d=2$ and $q=3$.\\
\begin{theorem}\label{thm-sierksma}
Sierksma's conjecture on the number of Tverberg partitions holds
for $q=3$ and $d=2$.
\end{theorem}

This paper is organized as follows: Section~\ref{sec-prel} comes with
a reminder of what is needed in the subsequent sections. In
Section~\ref{sec-prf-tvc}, we prove
Theorem~\ref{thm-constraint-graphs}. In Section~\ref{sec-conn-res},
we obtain the connectivity results for the chessboard-type complexes
needed in Section~\ref{sec-prf-tvc}. In Section~\ref{sec-number-tp},
we prove Theorems~\ref{thm-low-aff-improv} and~\ref{thm-sierksma}. 
\end{section}

\begin{section}{Preliminaries}\label{sec-prel}
Let's prepare our tools from topological combinatorics, and start with
some preliminaries to fix our notation, see also Matou\v{s}ek's 
textbook~\cite{matou03:_using_borsuk_ulam}. Let $k\geq -1$. A topological
space $X$ is {\it $k$-connected} if for every $l=-1,0,1,\ldots,k$,
each continuous map $f:S^l\ra X$ can be extended to a continuous map
$\bar{f}:B^{l+1}\ra X$.  Here $S^{-1}$ is interpreted as the empty set and
$B^0$ as a single point, so $(-1)$-connected means non-empty. We write
$\conn(X)$ for the maximal $k$ such that $X$ is $k$-connected.  There
is an inequality for the connectivity of the join $X*Y$ for
topological spaces $X$ and $Y$ which we use:
\begin{eqnarray}\label{conn}
\conn(X*Y)\geq \conn(X)+\conn(Y)+2;
\end{eqnarray}
see also~\cite[Section 4.4]{matou03:_using_borsuk_ulam}.

{\bf Deleted joins.} The {\it
  $n$-fold $n$-wise deleted join} of a topological space $X$ is
\[X^{*n}_{\Delta}:=X^{*n}\setminus\{{\textstyle\frac{1}{n}x_1\oplus\frac{1}{n}x_2
  \oplus \cdots\oplus\frac{1}{n}x_n\,}|\,\text{ $n$ of the $x_i\in X$ are
  equal}\}.\] We remove the diagonal elements from the
$n$-fold join $X^{*n}$.

For a simplicial complex $\Ka$ we define its {\it $n$-fold
  pairwise deleted join} as the following set of simplices:
\[\Ka^{*n}_{\Delta(2)}:=\{ F_1\uplus F_2\uplus\cdots\uplus F_n\in
\Ka^{*n}\,|\, F_1,F_2,\ldots,F_n \mbox{ pairwise disjoint}\}.\]

Both constructions show up in the proof of the topological Tverberg
theorem.
The $p$-fold pairwise deleted join of the $n$-simplex $\sigma^n$
is isomorphic to the $n+1$-fold join of a discrete space of $p$ points:
\begin{eqnarray}\label{eqn-iso-conf} 
(\sigma^n)^{*p}_{\Delta(2)}\cong ([p])^{*(n+1)}.
\end{eqnarray} 
In particular,
the simplicial complex $(\sigma^n)^{*p}_{\Delta(2)}$ is $n$-dimensional,
and $(n-1)$-connected.

\begin{figure}[h]
  \centering
  \includegraphics{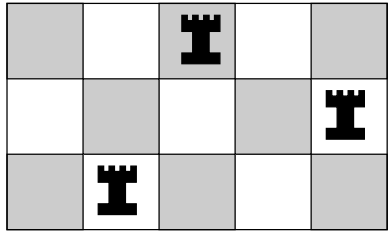}
  \caption{A maximal face of the chessboard complex $\Delta_{3,5}$.}
  \label{fig-chessboard}
\end{figure}
The {\it chessboard complex $\Delta_{m,n}$} is defined as the
simplicial complex $([n])^{*m}_{\Delta(2)}$. Its vertex set is the set
$[n]\times [m]$, and its simplices can be interpreted as placements of
rooks on an $n\times m$ chessboard such that no rook threatens any
other; see also Figure~\ref{fig-chessboard}. The roles of $m$ and $n$
are hence symmetric. $\Delta_{m,n}$ is an $(n-1)$-dimensional
simplicial complex with ${m\choose n}n!$ maximal faces for $m\geq n$.
See also Figure~\ref{fig-chessboard}, every maximal face corresponds
to a placement of $3$ rooks on a $3\times 5$ chessboard. Having 
equation~(\ref{eqn-iso-conf}) in mind, the chessboard complex $\Delta_{n,p}$
can be seen as a subcomplex of~$(\sigma^n)^{*p}_{\Delta(2)}$.

{\bf Nerve Theorem.} 
Another very useful tool in topological combinatorics is the nerve
theorem, e.~g.~it can be used to determine the connectivity of a given
topological space, or simplicial complex.  The {\it nerve $N(\Fe)$} of
a family of sets $\Fe$ is the abstract simplicial complex with vertex
set $\Fe$ whose simplices are all $\sigma\subset\Fe$ such that
$\bigcap_{F\in\sigma} F\not=\emptyset$. 

The nerve theorem was first obtained by Leray~\cite{leray45:_sur}, and
it has many versions; see Bj\"orner~\cite{bjoerner95:_topol} for a
survey on nerve theorems. \
\begin{theorem}[Nerve theorem]
 \label{thm-nerve-v-ii}
 For $k\geq 0$, let $\Fe$ be a finite family of subcomplexes of
 simplicial complex such that $\bigcap \Ge$ is empty or
 $(k-|\Ge|+1)$-connected for all non-empty subfamilies
 $\Ge\subset\Fe$.  Then the topological space $\|\bigcup\Fe\|$ is
 $k$-connected iff the nerve complex $\|N(\Fe)\|$ is $k$-connected.
\end{theorem}
Using Theorem~\ref{thm-nerve-v-ii} and induction, Bj\"orner,
Lov\'asz, Vre\'cica, and \v{Z}ivaljevi\'c proved in
\cite{bjoerner94:_chess} the following connectivity result for the
chessboard complex.
\begin{theorem}\label{thm-conn-chessboard}The chessboard complex
$\Delta_{m,n}$ is $(\nu-2)$-connected, for
\[ \nu :=\min\,\{ m,n,\lfloor \tfrac{1}{3}(m+n+1)\rfloor\} .\]
\end{theorem}

{\bf G-spaces and equivariant maps.} Let $(G,\cdot)$ be a
finite group with $|G|>1$. A topological space $X$ equipped with a
(left) $G$-action via a group homomorphism $\Phi:(G,\cdot)\ra
(\mbox{Homeo}(X),\circ)$ is a {\it $G$-space $(X,\,\Phi)$}.
Here Homeo$(X)$ is the group of
homeomorphisms on $X$, the product $\circ$ of two homeomorphisms $h_1$ and
$h_2$ is their composition. A continuous map $f$
between $G$-spaces $(X,\Phi)$ and $(Y,\Psi)$ that commutes with the
$G$-actions of $X$ and $Y$ is called a {\it $G$-map}, or an {\it
  equivariant map}. For  $x\in X$ the 
set $O_x=\{g\, x\,|\,g\in G\}$ is called the {\it
orbit} of $x$. A $G$-space $(X,\Phi)$ where every $O_x$ has at
least two elements is called {\it fixed point free}, i.~e.~no point
of X is fixed by all group elements.

The spaces $(\sigma^n)^{*q}_{\Delta(2)}$, $\Delta_{q,n}$, and
$(\R^n)^{*q}_{\Delta}$
are examples of $S_q$-spaces, where $S_q$ is the symmetric group
on $q$ elements. $S_q$ acts on all three spaces via permutation 
of the $q$ factors. For every subgroup $H$ of $S_q$, e.~g.~$\Z_q$,
or $(\Z_p)^r$ for prime powers $q=p^r$, 
an $S_q$-space is turned into a $H$-space
via restriction. In fact, $(\R^n)^{*q}_{\Delta}$ is a fixed point
free $(\Z_p)^r$-space for prime powers $q=p^r$, see for example
Hell~\cite[Lemma 5]{hell07:_tverb}.

It is one of the key steps in the equivariant method 
to prove that there is no $G$-map between
two given $G$-spaces. It is sufficient to prove that there is no
$H$-map between the $H$-spaces obtained via restriction, for a
subgroup $H$ of $G$. In the proof of the topological
Tverberg theorem for primes $q$ in the version 
of~\cite{matou03:_using_borsuk_ulam}, this is shown for
the subgroup $\Z_q$ via a $\Z_q$-index argument.

A less standard tool from equivariant topology is due to
Volovikov~\cite{volovikov96:_tverb_theor}. A cohomology $n$-sphere over
$\Z_p$ is a CW-complex having the same cohomology groups with
$\Z_p$-coefficients as the $n$-dimensional sphere~$S^n$. The space
$(\R^d)^{*q}_{\Delta}$ being homotopic to the $(d+1)(q-1)-1$-sphere 
is an example of a
cohomology $(d+1)(q-1)-1$-sphere over $\Z_p$, see for example 
Hell~\cite[Lemma 6]{hell07:_tverb}.
\begin{proposition}[Volovikov's Lemma]
  \label{lem-vol}Set $G=(\Z_p)^r$, and let $X$ and $Y$ be fixed point
  free $G$-spaces such that $Y$ is a finite-dimensional cohomology
  $n$-sphere over $\Z_p$ and $\tilde{H}^i(X,\Z_p)=0$ for all $i\leq
  n$. Then there is no $G$-map \mbox{from $X$ to $Y$}.
\end{proposition}
It is the key result in~\cite{volovikov96:_tverb_theor} to obtain 
Theorem~\ref{thm-ttt} for prime powers~$q$.\\

{\bf On Tverberg and Birch partitions.}  For 
Theorems~\ref{thm-low-aff-improv} and~\ref{thm-sierksma}, we have to
review some recent results for the affine setting of Tverberg's
theorem. A set of points in $\R^d$ is
{\it in general position} if the coordinates of all points are
independent over $\Q$. 
We have chosen this quite restrictive definition of general position 
for the sake of its brevity, see 
also~\cite{schoeneborn05:_topol_tverb_theor} for a less 
restrictive definition.
We need the following reformulation
of Lemma~2.7 from Sch\"oneborn and
Ziegler~\cite{schoeneborn05:_topol_tverb_theor}.
\begin{lemma}\label{obs-gen-pos-tp}
  Let $X$ be a set of $(d+1)(q-1)+1$ points in
  general position in $\R^d$. Then a Tverberg partition consists of:\\
  \textbullet\,\, Type~I: One vertex $v$, and $(q-1)$ many
  $d$-simplices containing $v$. \\
  \textbullet\,\, Type~II: $k$ intersecting simplices of dimension
  less than $d$, and $(q-k)$ $d$-simplices containing the
  intersection point for some $1<k\leq \min\{d,q\}$.
\end{lemma}
For $d=2$, a type~II partition consists of two intersecting segments,
and $q-2$ many triangles containing their intersection point. For
both types, the vertex resp.~the intersection point is a Tverberg
point.

Let $X$ be a set of $k(d+1)$ points in $\R^d$ for some $k\geq 1$. A
point $p\in\R^d$ is a {\it Birch point} of $X$ if there is a partition
of $X$ into $k$ subsets of size $d+1$, each containing $p$ in its
convex hull. The partition of $X$ is a {\it Birch partition for $p$}.
Let $B_p(X)$ be the number of Birch partitions of $X$ for $p$. If
$p$ is not in the convex hull of $X$, then clearly $B_p(X)=0$.

A Tverberg partition of a set of $(d+1)(q-1)+1$ points in $\R^d$ 
is an example of a Birch partition: For a type~I partition, one of the
points of this set is the Tverberg point. This point plays the role 
of the point $p$, and the remaining $(q-1)(d+1)$ points are partitioned
into $q-1$ subsets of size $d+1$. For a type~II partition, the intersection 
point is the Tverberg point which plays the role of the point $p$, and
the remaining points are again partitioned into subsets of size $d+1$.  
Now Theorem~\ref{thm-lower-aff-arb} follows from the following result
from Hell~\cite{hell07:_birch}.
\begin{theorem}\label{thm-number-birch-part}Let $d\geq 1$ and $k \geq 2$ be
  integers, and $X$ be a set of $k(d+1)$ points in $\R^d$ in general
  position with respect to the origin $0$. Then the following
  properties hold for $B_0(X)$:
\begin{enumerate}[\rm i)]
\item\label{it-birch-even}$B_0(X)$ is even.
\item\label{it-birch-lower}$B_0(X)>0\,\,\Longrightarrow\,\,B_0(X)\geq k!$
\end{enumerate}
\end{theorem}

\end{section}

\begin{section}{Proof of Theorem
    \ref{thm-constraint-graphs}}\label{sec-prf-tvc}

Figure~\ref{fig-constr-graphs5} shows all known elementary constraint
graphs for $q=5$, except for cycles on more than four vertices.  
In general, intersection graphs are disjoint unions
of elementary constraint graphs in the $1$-skeleton of $\sigma^{N}$.
For $q=2$, there are no constraint graphs. For $q=3$, a single edge
$K_2$ is the only elementary constraint graph.
\begin{figure}[h]
  \centering
  \includegraphics{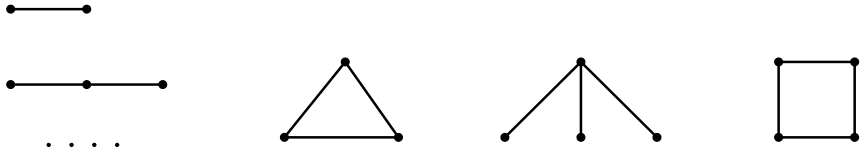}
  \caption{All known elementary constraint graphs for $q=5$.}
  \label{fig-constr-graphs5}
\end{figure}

\bprf (of Theorem~\ref{thm-constraint-graphs}) Set $N:=(d+1)(q-1)$,
and let $q>2$ be of the form $p^r$ for some prime number $p$. 
As in the proof of topological Tverberg theorem in the version
of~\cite{matou03:_using_borsuk_ulam}, we consider the
space $\Ka:=(\sigma^N)^{*q}_{\Delta(2)}$ as configuration
space. It models all possible partitions of the vertex set into
$q$ blocks: A maximal simplex of $\Ka$ encodes a (Tverberg) 
partition as shown in
Figure~\ref{fig-encoding}, and it can be represented as a hyperedge
using one point from each row of $\Ka$.

\begin{figure}[h]
  \centering
  \includegraphics{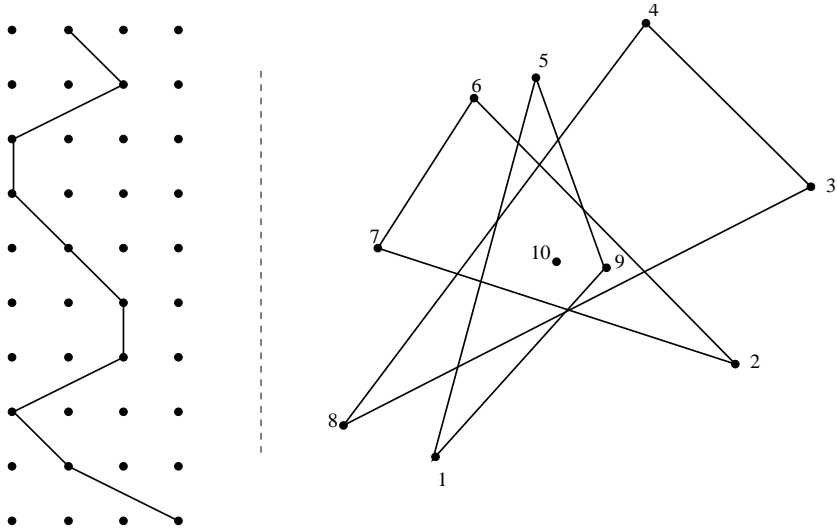}
  \caption{Maximal simplex of
    $(\sigma^N)^{*q}_{\Delta(2)}$ encoding a Tverberg partition.}
  \label{fig-encoding}
\end{figure}

Remember that $\|\Ka\|$ is $N-1$-connected. In the original proof of 
Theorem~\ref{thm-ttt},
the assumption that there is no Tverberg partition for $f$ leads to
the existence of a $(\Z_p)^r$-map $f^q:\|\Ka\| \ra (\R^d)^{*q}_{\Delta}$.
However, there is not such a map due to Volovikov's Lemma~\ref{lem-vol}. 
Hence a Tverberg partition exists for $f$.

In the following, we construct for each graph a good subcomplex $\La$ of
$\Ka$ such that: i) $\La$ is invariant
under the $(\Z_p)^r$-action, and ii) $\conn(\La)\geq N-1$.  Here {\it
  good} means that $\La$ does not contain any of Tverberg partitions
using an edge of our graph. As in the subsequent paragraph,
the assumption that
there is no Tverberg partition leads to
a $(\Z_p)^r$-map $f^q:\|\La\| \ra (\R^d)^{*q}_{\Delta}$. Finally
Volovikov's Lemma~\ref{lem-vol} implies a contradiction, and so
that there is a Tverberg partition not using any edge of our graph.
Hence, our
graph is a constraint graph.

Our construction of good subcomplexes is based in its
simplest case -- for $K_2$ -- on the following observation:
\begin{quote}If two points $i$ and $j$ end up in the same partition
  set, then the maximal face representing this partition uses one of
  the vertical edges between the corresponding rows $i$ and $j$ in
  $\Ka$.
\end{quote}
To prove the $K_2$ case, we have to come up with a subcomplex $\La$
that does not contain maximal simplices using vertical edges between
rows $i$ and $j$. Let $\La$ be the join of the chessboard complex
$\Delta_{2,q}$ on rows $i$ and $j$, and the remaining rows.
Figure~\ref{fig-encoding-constraint} shows this construction of $\La$
for $q=3$ and $d=2$. The chessboard complex $\Delta_{2,q}$ does not
contain any vertical edges. Moreover, $\La$ is $(\Z_p)^r$-invariant as
only the orbit of the vertical edges is missing. For the connectivity
of $\La$ see the next paragraph.

\ref{comp-const}) Construction of $\La$ for complete graphs $K_l$:
Let $i_1, i_2,\ldots, i_l$ be the corresponding rows of $\Ka$. $\La$
must not contain any maximal faces with vertical edges between any
two of these rows. The chessboard complex on these rows is such
a candidate. Let $\La$ be the join of the chessboard complex $\Delta_{l,q}$ on the corresponding $l$ rows, and the remaining rows:
\[ \La = \Delta_{l,q} * ([q])^{*(N+1-l)}.\] 
The subcomplex $\La$ is
closed under the $(\Z_p)^r$-action. Using Theorem
\ref{thm-conn-chessboard} on the connectivity of the chessboard
complex, and inequality (\ref{conn}) on the connectivity of the join, we
obtain:
\begin{eqnarray*}
\conn (\La)& \geq &\conn (\Delta_{l,q}) + \conn (([q])^{*(N+1-l)})+2\\
& \geq & \conn (\Delta_{l,q}) + N-l+1\\
& \geq & N-1.
\end{eqnarray*}
In the last step, we use that $\Delta_{l,q}$ is $(l-2)$-connected for
$2l<q+2$.

\begin{figure}[h]
  \centering
  \includegraphics{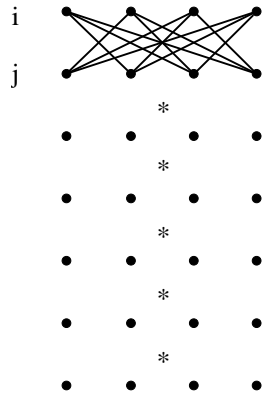}
  \caption{The construction of $\La$ for $K_2$.}
  \label{fig-encoding-constraint}
\end{figure}

\ref{bip-const}) Construction of $\La$ for complete bipartite graphs
$K_{1,l}$:
We first construct an $(\Z_p)^r$-invariant subcomplex
$C_{l,q}$ on the corresponding $l+1$ rows. For this, let $i$ be the
row that corresponds to the vertex of degree $l$, and $j_1, j_2,\ldots
j_l$ be the corresponding rows to the $l$ vertices of degree $1$. 
Let $C_{l,q}$ be the maximal induced subcomplex of $\Ka$ on the rows $i,
j_1, j_2,\ldots, j_l$ that does not contain any vertical edges
starting at a vertex of row $i$. Then $C_{l,q}$ is the union of $q$
many complexes $L_1, L_2,\ldots, L_q$, which are all of the form of
$\cone([q-1]^{*l})$.  Here the apex of $L_m$ is the $m$th vertex of
row $i$ for every $m=1, 2,\ldots, q$. In Figure~\ref{fig-complex-bip},
the maximal faces of the
complex $L_3$ are shown for $q=4$, and $l=2$.\\
Let $\La$ be the join of the complex $C_{l,q}$ and the remaining rows
of $\Ka$:
\[\La = C_{l,q}*([q])^{*(N-l)}.\]
Now $\La$ is good and $(\Z_p)^r$-invariant by construction.
Let's assume
\begin{eqnarray}\label{conn-constraint-complex}\conn( C_{l,q})\geq l-1
\end{eqnarray}
for $1<l<q-1$. The connectivity of $\La$ is then shown as above:
\begin{eqnarray*}
\conn (\La)& \geq &\conn (C_{l,q}) + \conn (([q])^{*(N-l)})+2\\
& \geq & \conn (C_{l,q}) + N-l\\
& \geq & N-1.
\end{eqnarray*}
We prove assumption (\ref{conn-constraint-complex}) in Lemma
\ref{lem-conn-constr-comp-bip} below.
\begin{figure}[h]
  \centering
  \includegraphics{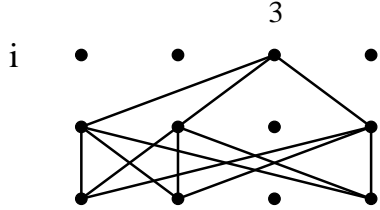}
  \caption{The complex $L_3$ for $q=4$ and $l=2$.}
  \label{fig-complex-bip}
\end{figure}

\ref{path-const}) Construction of $\La$ for paths $P_l$ on $l+1$
vertices: We construct recursively a good subcomplex $\La$ on
$l+1$ rows such that $\conn(\La)\geq l-1$. The case $l=1$ is covered
in the proof of~\ref{comp-const}) so that we can choose $\La$ to be
the complex $D_{2,q}:=\Delta_{2,q}$. For $l>1$, choose $\La$ to be the
complex $D_{l,q}$ which is obtained from $D_{l-1,q}$ in the following
way: Order the corresponding rows $i_1,i_2,\ldots,i_{l+1}$ in the
order they occur on the path. Take $D_{l-1,q}$ on the first $l$ rows.
A maximal face $F$ of $D_{l-1,q}$ uses a point in the last row $i_l$
in column $j$, for some $j\in[q]$. We want $D_{l,q}$ to be good so
that we cannot choose any vertical edges between row $i_l$ and
$i_{l+1}$.  Let $D_{l,q}$ be defined through its maximal faces: All
faces of the form $F\uplus\{k\}$ for $k\not=j$. Let $D_{l,q}^{k}$
be the subcomplex
 of all faces $D_{l,q}$ ending with $k$. Then
$D_{l,q}=\bigcup_{k=1}^qD_{l,q}^{k}$. In Figure~\ref{fig-good-pl} the
recursive definition of the complex $D_{l,5}^{2}$ is shown.
\begin{figure}[h]
  \psfrag{2}{$2$}
  \psfrag{a}{$D_{l-1,5}$}
  \centering
  \includegraphics{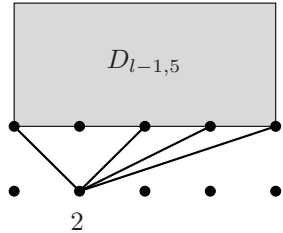}
  \caption{Recursive definition of $D_{l,5}^2$.}
  \label{fig-good-pl}
\end{figure}
The complex is $(\Z_p)^r$-invariant, and the connectivity
of $D_{l,q}$
\[ \conn(D_{l,q})\geq l-1\] is shown in Lemma
\ref{lem-conn-constr-comp-path} below using the decomposition
$\bigcup_{k=1}^qD_{l,q}^k$.

\ref{cycle-const}) Construction of $\La$ for cycles $C_l$ on $l$
vertices: Choose $\La$ to be the complex $E_{l,q}$ obtained from
$D_{l-1,q}$ on $l$ rows by removing all maximal simplices that use a
vertical edge between first and last row. The following result on the
connectivity of $E_{l,q}$ is shown in Lemma
\ref{lem-conn-constr-comp-cycle} below:
\[\conn(E_{l,q})\geq l-2.\]

\ref{union-const}) Construction of $\La$ for disjoint unions of
constraint graphs: For every graph component construct a complex on
the corresponding rows as above. Let $\La$ be the join of these
subcomplexes, and of the remaining rows.  Then $\La$ is a good
$(\Z_p)^r$-invariant subcomplex by the similar arguments as above. The
connectivity of $\La$ follows analogously from inequality (\ref{conn})
on the connectivity of the join. \eprf
\begin{remark*}

Figure~\ref{fig-k7-min-constr} comes with an example of a configuration 
of seven points in the plane showing that $P_2=K_{1,2}$ is not a constraint 
graph for $q=3$. This configuration is the outcome of a computer
program, see~\cite[Chapter 4]{hell06:_tverb_fract_helly} 
for details. The same program produced many planar point configurations
showing that $C_4$ is not a constraint graph for $q=4$.
\end{remark*}
\end{section}

\begin{section}{Connectivity for chessboard-type complexes}
\label{sec-conn-res}
The following three lemmas provide the connectivity results needed in
the proof of Theorem~\ref{thm-constraint-graphs}. Their proofs are
similar: Inductive on $l$, and Theorem~\ref{thm-nerve-v-ii}
is applied to the decompositions of the corresponding complexes that
were introduced in the proof of Theorem~\ref{thm-constraint-graphs}.

\begin{lemma}\label{lem-conn-constr-comp-bip}Let $q>2$, $d\geq 1$, and
  set $N=(d+1)(q-1)$. Let $C_{l,q}$ be the above defined subcomplex
  of $(\sigma^N)^{*q}_{\Delta(2)}$ for $1\leq l<q-1$. Then
  \[\conn(C_{l,q})\geq l-1.\]
\end{lemma}
\bprf In our proof, we use the decomposition of $C_{l,q}$ into
subcomplexes $L_1, L_2,\ldots$ $ L_q$ from above.

The nerve $\cal N$ of the family $L_1, L_2,\ldots, L_q$ is a
simplicial complex on the vertex set $[q]$. The intersection of $t$
many $L_{m_1}, L_{m_2},\ldots, L_{m_t}$ is $[q-t]^{*l}$ for $t>1$ so
that the nerve $\cal N$ is the boundary of the $(q-1)$-simplex.  Hence
$\cal N$ is $(q-3)$-connected.

Let's look at the connectivity of the non-empty intersections
$\bigcap_{j=1}^tL_{m_j}$. For $t=1$, every $L_m$ is contractible as it
is a cone. For $1<t<q-1$, the space $[q-t]^{*l}$ is $(l-2)$-connected,
and for $t=q-1$ the intersection is non-empty, hence its connectivity
is $-1$. All non-empty intersections $\bigcap_{j=1}^tL_{m_j}$ are thus
$(l-t)$-connected. The $(l-1)$-connectivity of $C_{l,q}$ immediately
follows from the nerve theorem using $q>2$, and
$l<q-1$. \eprf

\begin{lemma}\label{lem-conn-constr-comp-path}Let $q>3$, $d\geq 1$, and
  set $N=(d+1)(q-1)$.  Let $D_{l,q}$ be the above defined subcomplex
  of $(\sigma^N)^{*q}_{\Delta(2)}$ for $l\leq N$.  Then
  \[\conn(D_{l,q})\geq l-1.\]
\end{lemma}
\bprf In our proof, we use the decomposition of $D_{l,q}$ into
subcomplexes\linebreak $D^1_{l,q}, D_{l,q}^2,\ldots, D_{l,q}^q$ from
above.  We prove the following connectivity result by an induction on
$l\geq 1$:
\begin{eqnarray}\label{ineq-conn-dlq}
\conn(\bigcup_{j\in S} D_{l,q}^j)\geq l-1,\,\,\text{ for any
  $\emptyset\not=S\subset [q]$.}
\end{eqnarray} 

Let $l=1$, then $D_{1,q}=\bigcup_{j\in [q]}D_{1,q}^j$ is the
chessboard complex $\Delta_{2,q}$ which is $0$-connected for $q>2$.
The union of complexes $D_{1,q}^i$ is a union of contractible cones
which is $0$-connected. For $l\geq 2$, look at the intersection of $t>1$
many complexes $D_{l,q}^{i}$.  Let $T\subset [q]$ be the corresponding
index set of size $1<t<q-1$, and $\bar{T}$ its complement in $[q]$. Then
their intersections are
\begin{eqnarray}
\label{eq-intersec-dlq}
\bigcap_{j\in T}D_{l,q}^{j}&=&\bigcup_{j\in\bar{T}} D_{l-1,q}^j\,,\,\,\\
\label{eq-intersec-dlq-all-1}
\bigcap_{j\in [q]\setminus\{k\}}D_{l,q}^{j}&=& D_{l-1,q}^k
\cup D_{l-2,q}^k\,,\,\,\text{ and }\\
\label{eq-intersec-dlq-all}
\bigcap_{j\in [q]}D_{l,q}^{j}&=&\bigcup_{j\in [q]} D_{l-2,q}^j.
\end{eqnarray} 
The nerve $\cal N$ of the family $D_{l,q}^1, D_{l,q}^2,\ldots,
D_{l,q}^q$ is a simplicial complex on the vertex set $[q]$. The nerve
is the $(q-1)$-simplex, which is contractible.

For $l=2$, let's
apply the nerve theorem. For this, we have to
check that the non-empty intersection of any $t\geq 1$ complexes is
$(2-t)$-connected. Every $D_{2,q}^j$ is $1$-connected as it is a cone.
The intersection of $t=2$ many complexes is $0$-connected for $q>3$ by
equation~(\ref{eq-intersec-dlq}). Note that this is false for $q=3$.
The intersection of $t=3$ many complexes is non-empty. 

For $l=3$, we have to show that the non-empty intersection of 
any $t$ complexes is $(3-t)$-connected. Every $D_{3,q}^j$ 
is $2$-connected as it is a cone. The intersection of $t<q-1$ many 
complexes is $1$-connected by equation~(\ref{eq-intersec-dlq}).
The intersection of $t=q-1$ many complexes is a union of two cones
due to equation~(\ref{eq-intersec-dlq-all-1}). The intersection of
these two cones is:
\[ D_{2,q}^k\cap D_{1,q}^k=[q]\setminus\{k\},\]
which is non-empty. Using the nerve theorem, we obtain for their union:
\[\conn(D_{2,q}^k\cup D_{1,q}^k)\geq 0\geq 3-(q-1)\,\,\text{ for }q\geq 4.\]
The intersection of $t=q\geq 4$ many complexes is non-empty by 
equation~(\ref{eq-intersec-dlq-all}).

Let now $l>3$, we apply again the nerve theorem to obtain inequality
(\ref{ineq-conn-dlq}). It remains to check that the non-empty
intersection of any $t$ complexes is $(l-t)$-connected. The
complex $D_{l,q}^j$ is $(l-1)$-connected as it is a cone for every
$j\in [q]$.  The intersection of any $1<t<q-1$ complexes is
$(l-2)$-connected by equation (\ref{eq-intersec-dlq}) and by assumption.
The intersection of $t=q-1$ many complexes is a union of two cones
due to equation~(\ref{eq-intersec-dlq-all-1}). The intersection of
these two cones is:
\[ D_{l-1,q}^k\cap D_{l-2,q}^k=\bigcup_{j\in [q]\setminus\{k\}}D_{l-3,q}^{j},\]
which is $(l-4)$-connected by assumption. 
Using the nerve theorem, we obtain for their union:
\[\conn(D_{l-1,q}^k\cup D_{l-2,q}^k)\geq l-3\geq l-(q-1)\,\,\text{ for }q\geq 4.\]
The intersection of $q$ many complexes is $(l-3)$-connected by 
equation~(\ref{eq-intersec-dlq-all}) and by assumption. \eprf

\begin{lemma}\label{lem-conn-constr-comp-cycle}Let $q>4$, $d\geq 1$, and
  set $N=(d+1)(q-1)$.  Let $E_{l,q}$ be the above defined subcomplex
  of $(\sigma^N)^{*q}_{\Delta(2)}$ for $l\leq N+1$.  Then
  \[\conn(E_{l,q})\geq l-2.\]
\end{lemma}
\bprf The proof is similar to the proof of Lemma
\ref{lem-conn-constr-comp-path}. The case $l=3$ has
already been settled in the proof of case~\ref{comp-const}) of Theorem
\ref{thm-constraint-graphs}. The cases $l=4, 5$ are analogous for $q\geq 5$, 
but need some tedious calculations. Observe that the inductive argument
in the proof of Lemma~\ref{lem-conn-constr-comp-path} also works for
$E_{l,q}$, which was obtained from $D_{l-1,q}$ by removing some
maximal faces.

Let's describe the differences to the proof of 
Lemma~\ref{lem-conn-constr-comp-path}. We consider the decomposition 
$E^1_{l,q}, E^2_{l,q},\ldots ,E^q_{l,q}$ of $E_{l,q}$. Here $E^i_{l,q}$
is the complex that is obtained from $D^i_{l-1,q}$ by removing all maximal
faces that contain the $i$th vertex of the first row. In 
Figure~\ref{fig-good-cycle} the complex $E^1_{l,5}$ is shown: Any 
face of $D^1_{l-1,q}$ containing one of the broken edges is removed. 
\begin{figure}[h]
  \psfrag{1}{$1$}
  \psfrag{a}{$D_{l-3,5}$}
  \centering
  \includegraphics{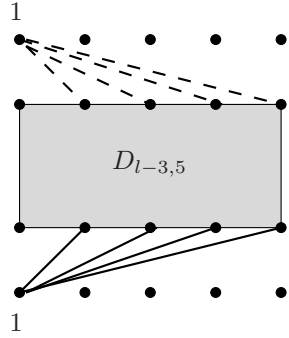}
  \caption{Subcomplex $E_{l,5}^1$ of $D_{l-1,5}^1$.}
  \label{fig-good-cycle}
\end{figure}

The intersection of this family is non-empty, in fact:
\begin{eqnarray}
\label{eq-intersec-e-q}
\bigcap_{i=1}^q E^i_{l,q}= D_{l-4,q}\,\,\text{ for }q\geq 5.
\end{eqnarray}
Thus its nerve is a simplex. Using the nerve theorem it remains to show
that the intersection of $t\geq 1$ complexes is $(l-2-t+1)$-connected. For $t=1$, the complex $E^i_{l,q}$ is a cone. For $t=q$,
this follows from equation~(\ref{eq-intersec-e-q}). For $1<t<q$, this follows
as in the proof of Lemma~\ref{lem-conn-constr-comp-path} from the equations:
\begin{eqnarray}
\label{eq-intersec-e-q-1}
\bigcap_{i\in [q]\setminus\{k\}} E^i_{l,q}= 
\tilde{D}^{k,[q]\setminus\{k\}}_{l-2,q}\cup \tilde{D}^{k,[q]\setminus\{k\}}_{l-3,q} \,\,\text{, and}\\ 
\label{eq-intersec-e-i}
\bigcap_{i\in T}^q E^i_{l,q}= \bigcup_{{i\in\bar{T}}}\tilde{D}^{i,T}_{l-2,q}\,\, 
\text{ for $T\subset [q]$ and $1<|T|<q-1$},
\end{eqnarray}
where $\tilde{D}^{i,S}_{l,q}$ is the following subcomplex of $D^i_{l,q}$ 
for $S\subset [q]$: Delete all faces that contain a vertex in $S$ of the 
first row. In other words $\tilde{D}^{i,\{i\}}_{l,q}=E^i_{l+1,q}$, see also
Figure~\ref{fig-e-i} for equation~(\ref{eq-intersec-e-i}). 
There any face containing a broken edge is deleted from $D^i_{l,q}$.
\begin{figure}[h]
  \psfrag{a}{$D_{l-4,5}$}
  \centering
  \includegraphics{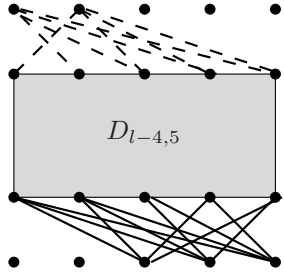}
  \caption{Equation~(\ref{eq-intersec-e-i}): $\bigcap_{i\in \{1,2\}}^q E^i_{l,5}= \bigcup_{{i\in\{3,4,5\}}}\tilde{D}^{i,\{1,2\}}_{l-2,5}$}
  \label{fig-e-i}
\end{figure}


Using again the nerve theorem, one then shows the necessary connectivity 
results for equations~(\ref{eq-intersec-e-q-1}) and~(\ref{eq-intersec-e-i}). 
This can be done for $q\geq 5$, inductively on $l\geq 5$:
\[ \conn(\tilde{D}^{k,[q]\setminus\{k\}}_{l-2,q}\cup \tilde{D}^{k,[q]\setminus\{k\}}_{l-3,q})\geq l-4, \] 
and for $T\subset [q]$, $1< |T| < q-1$:
\[\conn(\bigcup_{{i\in\bar{T}}}\tilde{D}^{i,T}_{l-2,q})\geq l-3,\text{ and }\,\,\conn(\bigcup_{{i\in T}}\tilde{D}^{i,T}_{l-2,q}))\geq l-3.\]
\eprf
\end{section}

\begin{section}{On the number of Tverberg partitions}\label{sec-number-tp}
  In this section, we start with the proof of 
  Theorem~\ref{thm-low-aff-improv}.  In the proof we apply
  Theorem~\ref{thm-constraint-graphs} on Tverberg partitions with
  constraints.  Using a similar approach, we then settle
  Sierksma's conjecture for $d=2$ and $q=3$.\\

  Having Theorem~\ref{thm-number-birch-part} in mind, we rise the
  following question:
\begin{quote}Is there a non-trivial lower bound for the number of \linebreak 
Tverberg points?
\end{quote}
In general, the answer is NO. Sierksma's well--known point
configuration has exactly one Tverberg point which is of type~I.  This
together with Theorem~\ref{thm-number-birch-part} 
leads to the term $(q-1)!$ in the lower bound of Theorem
\ref{thm-low-aff-improv}. But under the assumption that there are no
Tverberg points of type~I, we obtain a non-trivial lower bound for the
number of Tverberg points. The constant $c_{d,q}$ is in fact a lower
bound for the number of Tverberg points, assuming that there is none of
type~I. The factor $(q-d)!$ is due to the fact that we cannot predict
what kind of type~II partition shows up.

\bprf (of Theorem~\ref{thm-low-aff-improv}) Let $X$ be a set of
$(d+1)(q-1)+1$ points in $\R^d$, and $p_1$ is a Tverberg point which
is not of type~I.  The Tverberg point $p_1$ is the intersection point
of $\bigcap_{i=1}^k\conv(F_i^1)$, where $k\in\{2,3,\ldots,d\}$.
Choose an edge $e_1$ in some $F_i$, and apply Theorem
\ref{thm-constraint-graphs} with constraint graph $G_1=\{e_1\}$.  Then
there is a Tverberg partition that does not use the edge $e_1$ so that
there has to be second Tverberg point $p_2$.  Now add another edge
$e_2$ from the corresponding $F_i^2$ to the constraint graph $G_1$,
and apply again Theorem~\ref{thm-constraint-graphs} with constraint
graph $G_2=\{e_1, e_2\}$.  Hence there is another Tverberg point $p_3$
and so on. This procedure depends on the choices of the edges, and
whether $G_i$ is still a constraint graph.\\
Figure~\ref{fig-k7-min-constr} shows an example for $d=2$ and $q=3$: A
set of seven points in $\R^2$. There are exactly four Tverberg points
-- highlighted by small circles -- in this example. A constraint graph
-- drawn in broken lines -- can remove only three among them.
\begin{figure}[h]
  \centering
  \includegraphics{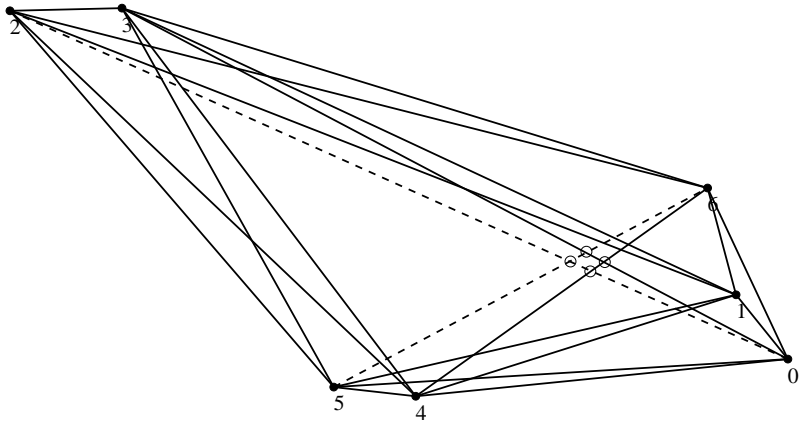}
  \caption{A set of $7$ points in the plane together with a constraint graph.}
  \label{fig-k7-min-constr}
\end{figure}\\
Constraint graphs for $q$ are also constraint graphs for the subsequent
prime power~$q'$ so that
our constant $c_{d,q}$ is weakly increasing in $q$.  The constant
$c_{d,q}$ also depends on $d$ as the simplex $\sigma^{(d+1)(q-1)}$
grows in $d$. 

It remains to prove $c_{2,3}>3$. For this, suppose we have 
three Tverberg partitions of type~II for the set 
$\{a,b,c,d,e,f,g\}$ of seven points in $\R^2$. 

If some edge, e.~g.~$\{a,b\}$, belongs to two partitions, we could
find an edge in the third partition disjoint with $\{a,b\}$. The
union of these two edges is a constraint graph.

If no edge belongs to two partitions, we have up to permutation
the Tverberg partitions $\{a,b,c\},\{d,e\},\{f,g\}$ and 
$\{a,d,f\},\{b,e\},\{c,g\}$and the third partition could be either
$\{a,e,g\},\{b,d\},\{c,f\}$ or $\{b,d,g\},\{a,e\},\{c,f\}$. 
In the former case the constraint graph $\{b,c\},\{d,f\},\{e,g\}$
contains an edge from
every partition, and shows that there has to be a fourth Tverberg
partition. In the later case, the same is true for the graph 
$\{b,c\},\{a,f\},\{d,g\}$. 
\eprf

Up to now, we have not been able to determine the exact value of
$c_{d,q}$ for $d>2$ or $q>3$, as there are just too many
configurations to look at. A similar -- in general smaller -- 
constant exists in the setting of the topological Tverberg
theorem.\\

{\bf On Sierksma's conjecture.} For $d=2$ and $q=3$,
Theorem~\ref{thm-low-aff-improv} settles Sierksma's conjecture for
sets having no type~I partition.  $c_{2,3}=4=((q-1)!)^d$ implies that
there are at least four different Tverberg partitions.  It remains to show
Sierksma's conjecture for planar set of seven points having
i)~only type~I partitions, and ii)~for sets with both partition types.

\bprf (of Theorem~\ref{thm-sierksma})
Case~i). There is at least one Tverberg point coming
with two partitions due to Theorem~\ref{thm-number-birch-part}.  It
remains to show that there is one more Tverberg partition, as evenness
implies the existence of the missing fourth one. Let $v$ be the
Tverberg point so that $\{v\},\{a,b,c\},\{d,e,f\}$ forms one of the
two Tverberg partitions.  Then the other Tverberg partition is of the
form $\{v\},\{a,b,d\},\{c,e,f\}$. Choosing for example
the edge $\{a,b\}$ as constraint graph completes our proof.
This is not the only possible choice for $G$.

Case~ii). There is again at least one Tverberg point
$v$ coming with two partitions of type~I: $\{v\},\{a,b,c\},\{d,e,f\}$
and $\{v\},\{a,b,d\},\{c,e,f\}$. The edge $\{a,b\}$ belongs to both
of these partitions. In the third partition of type~II, the points
$a$ and $b$ could belong to two sets of the partition. Choose any 
edge from the third set of this partition. It is disjoint with the
edge $\{a,b\}$, and together with it forms the constraint graph
showing that there has to be a fourth Tverberg partition.
\eprf
\end{section}

\section*{Final remarks}
Let's end with a list of problems on possible extensions of our
results. The first problem aims in the direction of finding similar
good subcomplexes. The second problem asks whether it is possible to
show the Tverberg theorem with constraints for affine maps,
independent of the fact that $q$ is a prime power. Moreover, we
conjecture that this method can be adapted to the setting of the
colorful Tverberg theorem.
\begin{problem*}\label{prob-constraint-graphs} Determine the
  class ${\cal CG}_{q,d}$ of constraint graphs. Find graphs that are
  not constraint graphs. Which of the constraint graphs are maximal?\\
  Show that cycles $C_l$ are constraint graphs for $q=4$, and $l\geq 5$.
\end{problem*}
\begin{problem*}\label{prob-cg-arbitrary} Identify constraint graphs for
 arbitrary $q\geq 2$, especially for affine maps.
\end{problem*}
\begin{problem*}\label{prob-lower-ctp} Find  good subcomplexes
  in the configuration space\linebreak $(\Delta_{2q-1,q})^{*d+1}$ of
  the colored Tverberg theorem to obtain a lower bound for the number
  of colored Tverberg partitions, and a colored Tverberg theorem with
  constraints.
\end{problem*}
Here a {\it good} subcomplex $(\Delta_{2q-1,q})^{*d+1}$ is again
$(\Z_p)^r$-invariant, and at least $((d+1)(q-1)-1)$-connected.
Constructing good subcomplexes in this setting requires more care than
for the topological Tverberg theorem. One possibility to construct
good subcomplexes is to identify $d+1$ many $(\Z_p)^r$-invariant
subcomplexes $\La_i$ in the chessboard complex $\Delta_{2q-1,q}$ such
that
\[ \sum_{i=1}^{d+1}\conn(\La_i)\geq (d+1)(q-3)+1. \] The join of the
$\La_i$'s is then a good subcomplex in $(\Delta_{2q-1,q})^{*d+1}$.
Looking at the proof for the connectivity of the chessboard complex,
and studying $\Delta_{2q-1,q}$ for small $q$ via the mathematical
software system polymake~\cite{joswig05:_geometric}, suggests that one
obtains subcomplexes $\La_i$ by removing a non-trivial number of
orbits of maximal faces. 

The last problem was suggested to me by G\'abor Simonyi.
\begin{problem*}Identify constraint hypergraphs.
\end{problem*}
Here a constraint hyperedge is a set of at least 3 vertices. All
vertices can not end up in the same block, but any subset can.
Forbidding a hyperedge of $n$ vertices is therefore weaker than
forbidding a complete graph $K_n$.\\

{\bf Acknowledgments.} The results of this paper are part of my PhD
thesis~\cite{hell06:_tverb_fract_helly}. I would like to thank
Juliette Hell, G\"unter M.~Ziegler, and Rade \v{Z}ivaljevi\'c for
many helpful discussions. Let me also thank the referees
for their insightful comments and corrections, which led 
to a substantial improvement of the paper.



\end{document}